\documentclass[11pt]{article}

\usepackage{amsfonts}

\newtheorem{Theorem}{Theorem}[section]
\newtheorem{Lemma}[Theorem]{Lemma}
\newtheorem{Proposition}[Theorem]{Proposition}
\newtheorem{Corollary}[Theorem]{Corollary}
\newtheorem{Def1}{Definition}[section]
\newtheorem{Exe1}{Example}[section]

\newcommand{\qed}{\hfill $\Box$ \hfill \\}
\long\def\symbolfootnote[#1]#2{\begingroup%
\def\thefootnote{\fnsymbol{footnote}}\footnote[#1]{#2}\endgroup}

\title{  \sc Braid action on derived category of Nakayama algebras}
\author{ Intan Muchtadi-Alamsyah}
\date{}

\begin{document}
\maketitle

\begin{abstract}
\noindent
We construct an action of a braid group associated to a complete graph on the derived category of a certain symmetric Nakayama algebra which is also a Brauer star algebra with no exceptional vertex. We connect this action with the affine braid group action on Brauer star algebras defined by Schaps and Zakay-Illouz. We show that for Brauer star algebras with no exceptional vertex, the action is faithful.
\end{abstract}

\section*{Introduction}
\symbolfootnote[0]{ 2000 {\it Mathematics Subject Classification} 16D90,16G20,18E30,20F36\\The author was supported by EPSRC Grant number GR/S35387/01.}

For a commutative ring $R$ we shall study the group of self-equivalences of the derived category of an $R$-algebra $A$ which is projective as an $R$-module. We denote such a group by $TrPic(A)$ and call it the {\it Derived Picard group.}

Rouquier and Zimmermann \cite{RouZi}, motivated by Brou\'e's conjecture, studied the $TrPic(A)$ of Brauer tree algebras $A$ (with no exceptional vertex) by defining an action of the Artin braid group $B(A_n)$ on $D^b(A)$, the bounded derived category of $A.$ Khovanov and Seidel \cite{Khovanov-Seidel}, motivated by mirror symmetry, defined an action of $B(A_n)$ on bounded derived category of certain algebras similar to  Brauer tree algebras, which turn out to be related to certain classes of representations of simple Lie algebras, and they proved that the action is faithful.

In this paper we will define an action of $B(K_n),$ the braid group
whose associated Coxeter group is given by the complete graph on $n$
vertices, on $D^b(N^n_n)$, the bounded derived category of the
Nakayama algebra $N^n_n$ which is a basic Brauer star algebra with
no exceptional vertex.

We will also  define a group homomorphism $\eta$ from $B(K_n)$ to $B(A_n)$ and give a categorification of $\eta.$ This categorification provides a link between the action of $B(K_n)$ on $D^b(N_n^n)$ and the faithful action of $B(A_n)$ on derived category of Brauer tree algebras defined by Rouquier and Zimmermann \cite{RouZi} and Khovanov and Seidel \cite{Khovanov-Seidel}.

In \cite{Schaps-Zakay-Illouz} Schaps and Zakay-Illouz defined an action of affine braid group $B(\tilde{A}_{n-1})$ on derived category of Brauer star algebras and left open the question whether this action is faithful.

We will provide a positive answer to this question for Brauer star algebras with no exceptional vertex. We will use the categorification of $\eta$ stated above, another categorification of braid group embeddings defined in  \cite{Kent-Peifer}, \cite{Graham-Lehrer} and \cite{Crisp} and  the faithfulness of the action in \cite{Khovanov-Seidel} and \cite{RouZi}.

The paper is organized as follows. \\
In Section 1, we review Rickard's Morita theory for derived equivalence and we give the definition of derived Picard groups $TrPic$ and of Nakayama algebras.

In Section 2, we recall the Nakayama algebra $N^n_n$ which is also a basic Brauer star algebra with $n$ edges and no exceptional vertex. We define explicit elements of $TrPic(N^n_n).$

The action of the braid group $B(K_n)$ on $D^b(N^n_n)$ is defined in
Section 3. And in Section 4 we show that this action is not faithful
using a categorification of the group homomorphism $\eta.$

In the last section, we use the results of sections 2, 3 and 4 to show that the action defined in \cite{Schaps-Zakay-Illouz} is faithful for Brauer star algebras with no exceptional vertex. \\

{\bf Acknowledgement} I would like to thank Robert Marsh and Steffen K\"onig for numerous discussions on studying the braid action on derived categories. I also thank Alexander Zimmermann for encouraging me to prove the faithfulness of the action defined in \cite{Schaps-Zakay-Illouz}.\\

\section{Derived Picard groups and Nakayama algebras}

\subsection{Derived equivalence}

Let $A$ and $B$ be algebras over a commutative ring $R$ and assume that $A$ and $B$ are projective as modules over $R.$ The bounded derived category $D^b(A)$ is the category  whose objects are complexes of finitely generated projective modules which are bounded to the right and which have nonzero homology only in finitely many degrees. Morphisms are morphisms of complexes up to homotopy. We refer to \cite{KonigZimmermann1685} for the definition of derived categories.

Denote by $A-per$  the full subcategory of $D^b(A)$  consisting of the {\it perfect complexes}, i.e. bounded complexes of finitely generated projective $A$-modules.

Rickard \cite{Rickarddereqderfun} and Keller \cite{Kellerremark} have given a necessary and sufficient criterion for the existence of derived equivalences between two rings $A$ and $B.$

\begin{Theorem}(Rickard \cite{Rickarddereqderfun} and Keller \cite{Kellerremark}) \\ The following statements are equivalent :
\begin{enumerate}
\item $D^b(A) \simeq D^b(B)$ as triangulated categories.
\item There is a complex $T$ in $A-per$  such that
\begin{enumerate}
\item $Hom(T,T[i]) = 0$ for $i \neq 0,$
\item $A-per$ is generated by $T,$ and
\item $End_{D^b(A)}(T) \simeq B.$
\end{enumerate}
\item There is a bounded complex $X$ in $D^b(A \otimes_R B^{op})$ whose restrictions to $A$ and to $B^{op}$ are perfect and a bounded complex $Y$ in $D^b(B \otimes_R A^{op})$ whose restrictions to $B$ and $A^{op}$ are perfect such that
\[ X \otimes^L_B Y \simeq A \;\; \hbox{in} \;\; D^b(A \otimes A^{op}) \;\; \hbox{and} \;\; Y \otimes^L_A X \simeq B \;\; \hbox{in} \;\; D^b(B \otimes B^{op}). \]
\end{enumerate} \end{Theorem}
A complex $T$ satisfying the condition in 2 is called a {\it tilting complex} for $A$ and the complexes $X$ and $Y$ in 3 are called  {\it two-sided tilting complexes} inverse to each other. The image and the pre-image of an  indecomposable projective $B$-module under derived equivalence $D^b(A) \simeq D^b(B)$ is a {\it partial tilting complex} in $A$, i.e. a complex satisfying the condition in 2(a). \\

We define the {\it Derived Picard group} of $A$ as
\[\begin{array}{ccc}
TrPic(A)&  = & \left\{ \begin{array}{c}  \hbox{isomorphism classes of two-sided tilting complexes} \\ \hbox{ in  } D^b(A \otimes_R A^{op}) \end{array} \right\} \end{array}\] where the product of the classes of $X$ and $Y$ is given by the class of $X \otimes_A Y.$

It is clear that if  $A$ and $B$  are projective as $R$-modules then $D^b(A) \simeq D^b(B)$ implies \[TrPic(A) \simeq TrPic(B).\]

 The Picard group $Pic(A)$ (i.e. the group of isomorphism classes of invertible $A \otimes A^{op}$-modules $M$ where $M \otimes_A M^* \simeq A$)  can be embedded into $TrPic(A)$ by sending an invertible $A \otimes A^{op}$-module $M$ to a stalk complex concentrated in degree 0.

\subsection{Nakayama algebras}

 The self-injective Nakayama algebra $N^n_m$ is the path algebra over a field $K$ of the following quiver
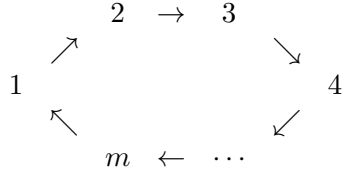
\[\begin{array}{ccccccc} &           & 2 & \rightarrow & 3      &          & \\
                         & \nearrow  &   &             &        & \searrow & \\
                       1 &           &   &             &        &          & 4 \\
                         & \nwarrow  &   &             &        & \swarrow &  \\
                         &           & m & \leftarrow  & \cdots    & & \end{array} \]
modulo the ideal generated by all compositions of $n+1$ consecutive arrows. By \cite{Asashibadereqclass} $N^n_m$ is a representative of the class of standard self-injective algebras of finite representation type associated to the Dynkin diagram of type $A_n.$

The algebra $N^n_m$ is symmetric if and only if $m$ divides $n.$ In this case the Nakayama algebra $N^n_m$ is a basic Brauer star algebra with $n$ edges and multiplicity $n/m$, a particular Brauer tree algebra where all the edges are adjacent to the exceptional vertex. A Brauer tree algebra with $n$ edges and multiplicity $n/m$ is (up to Morita equivalence) a $p$-block with cyclic defect group of order $p^d = m + 1$ where the number of its  simple modules is $n.$ (See \cite{Alperin} for the definition of Brauer tree algebras).

By \cite[Theorem 4.2]{Rickarddercatstabeq}, up to derived equivalence, a Brauer tree algebra is determined by the number of edges of the Brauer tree and the multiplicity of the exceptional vertex. Hence, an arbitrary Brauer tree algebra  is derived equivalent to a Brauer star algebra associated with a star having the same number of edges.

If $m = n,$ the Nakayama algebra $N^n_n$ is the trivial extension algebra of the path algebra of the quiver
\[ 1 \rightarrow 2 \rightarrow 3 \rightarrow \cdots \rightarrow n-1 \rightarrow n. \]
It is a Brauer star algebra with $n$ edges and multiplicity 1. \\

{\bf Notation} : From now on, $A$ is the Nakayama algebra $N^n_n$ over $K$ and $B$ is the Brauer tree algebra associated to a line without exceptional vertex. For projective $A$-modules we use the letter $P$ and for projective $B$-modules we use the letter $Q.$ We use the letters $F,$ $R,$ and $H$ to denote two-sided tilting complexes (elements of $TrPic$) and the letters $\mathcal{F},$ $\mathcal{R},$ and $\mathcal{H}$ for the corresponding auto-equivalences. The symbol $\otimes$ means $\otimes_K.$

\section{Two-sided tilting complexes for Nakayama algebra $N_n^n$}

Denote by $(i_1 \;\; i_2 \;\; \cdots i_{l+1})$ the path starting at $i_1$ and ending at $i_{l+1}.$ Note that we have $i_{k+1} - i_k = 1 \;(mod \;n)$ for all $1 \le k \le l.$

Let $A$ be the Nakayama algebra $N_n^n.$ The algebra $A$ is the path algebra of the quiver
  \[\begin{array}{ccccccc} &           & 2 & \rightarrow & 3      &          & \\
                         & \nearrow  &   &             &        & \searrow & \\
                       1 &           &   &             &        &          & 4 \\
                         & \nwarrow  &   &             &        & \swarrow &  \\
                         &           & n & \leftarrow  & \cdots    & & \end{array} \]
modulo the relations :
$(i \;\; i+1 \;\; ... \;\; n \;\; 1 \;\; 2 \;\; ... \;\; i \;\; i+1) = 0$ for all $1 \leq i \leq n.$ \\

The paths $(i)$ of length 0 are mutually orthogonal idempotents and their sum is the unit element.
The Loewy series of the  indecomposable projective modules $P_i = A(i)$ are as follows : {\small
\[ P_i = \begin{array}{c} S_i \\ S_{i+1} \\ \vdots \\ S_n \\ S_1 \\ S_2 \\ \vdots \\ S_i \end{array} \]}

The dimension of homomorphisms between projective modules is given by :
\[ dim_K Hom_A(P_i,P_j) = \left\{ \begin{array}{cc} 2 & \hbox{if } i = j, \\ 1 & \hbox{if } i \neq j. \end{array} \right.  \]
Define $_iP = (i)A \cong Hom_A(P_i,A).$ Define the complexes \[ F_i :=  0 \rightarrow P_i \otimes \; _iP \stackrel{\alpha_i}{\longrightarrow} A \rightarrow 0 \] and \[F_i' := 0 \rightarrow A \stackrel{\beta_i}{\longrightarrow} P_i \otimes \; _iP  \rightarrow 0 \] where $A$ is in degree 0 and
\[ \begin{array}{ccccccc} \alpha_i & : & P_i & \otimes & _iP & \rightarrow & A \\
                                   &   & (i) & \otimes & (i) & \mapsto     & (i) \end{array} \]
\begin{eqnarray*} \beta_i  :  & A \rightarrow & P_i \otimes \; _iP \\
                                   & 1  \mapsto     & (i \; \cdots \; n \; 1 \; \cdots \; i) \otimes (i) \; + \;(i) \otimes (i \; \cdots \; n \; 1 \; \cdots \;i) \\ &  & + \;(i+1 \; \cdots \; n \; 1 \; \cdots \; i) \otimes (i \; i+1) \\ & & + \; (i+2 \; \cdots \; n \; 1 \; \cdots \; i) \otimes (i \; i+1 \; i+2) \\ & & + \; \cdots \; + \; (n \; 1 \; \cdots \; i) \otimes (i \; i+1 \; \cdots \; n) \\ &  & + \;(1 \; \cdots \; i) \otimes (i \; \cdots \;n \;1)\; + \; (2 \; \cdots \; i) \otimes (i \; \cdots \; n \; 1 \; 2) \; \\ &  & + \; \cdots \; + \;(i-1 \; i) \otimes (i \; \cdots \;n \; 1 \; \cdots  \;i-1). \end{eqnarray*}

\begin{Theorem} The complexes $F_i$ and $F_i'$ are elements of $TrPic(A).$ \end{Theorem}

\noindent {\bf Proof} We will show that $F_i$ and $F_i'$ are two-sided tilting complexes and inverses to each other and to do that, we will show that $F_i \otimes_A F_i'$ is homotopy equivalent to $A$ as complex of $A-A-$bimodules. (See \cite[Theorem 6]{Rouquier}). We follow the method of  \cite[Theorem 4.1]{RouZi} and \cite[Proposition 2.4]{Khovanov-Seidel}. We have
\[ F_i \otimes_A F_i' = (0 \rightarrow P_i \otimes \; _iP \stackrel{d^1}{\longrightarrow} A \oplus (P_i \otimes U \otimes \; _iP) \stackrel{d^0}{\longrightarrow} P_i \otimes_A \; _iP \rightarrow 0) \] where $U$ is the 2-dimensional space $\;_iP \otimes_A P_i$ with a basis $u_1 = (i) \otimes (i)$ and $u_2 = (i \cdots i) \otimes (i),$ and the differentials are given by
\[ d^1 = \alpha_i + \tau \;\; \hbox{and} \;\;  d^0 = (\beta_i, -\delta) \; \hbox{where}\]
\begin{eqnarray*} \tau(x \otimes y) & = & x \otimes u_1 \otimes (i \cdots i)y + x \otimes u_2 \otimes y ,\\
\delta(x \otimes u_1 \otimes y) & = & x \otimes y ,\\
\delta(x \otimes u_2 \otimes y) & = & x(i \cdots i) \otimes y. \end{eqnarray*}

The map $d^0$ is surjective since $\delta$ is surjective. Therefore, since $P_i \otimes \; _iP$ is projective, $d^0$ is a split surjection.

Denote by $\mathcal{B}_i$ (resp. $_i\mathcal{B}$) the basis of $P_i$ (resp. of $_iP$) where  \[ \mathcal{B}_i = \{ (i), (i-1\; i), (i-2\;i-i\;i), \cdots, (n \;1\;2 \cdots i), \cdots, (i \; i+1 \; \cdots n \;1 \cdots i)\} \; \hbox{and}\]  \[ _i\mathcal{B} = \{(i), (i \;i+1), (i\;i+1, i+2), \cdots, (i \; i+1 \cdots n \; 1) , \cdots, (i \; i+1 \cdots n \; 1 \cdots i)\}.\]
Let $x \in \mathcal{B}_i$ and $y \in \;_i\mathcal{B}.$ Then,
\begin{eqnarray*} \tau(x \otimes y) & = & \left\{ \begin{array}{ccc} x \otimes u_1 \otimes (i \cdots i) + x \otimes u_2 \otimes (i) & \hbox{if}& y = (i), \\ x \otimes u_2 \otimes y & \hbox{if} & y \neq (i). \end{array} \right. \end{eqnarray*}

Therefore, if $x \otimes y \in P_i \otimes \;_iP,$ then $x \otimes y \neq 0$ implies $\tau(x \otimes y) \neq 0.$ This means that $\tau$ is injective, therefore so is $d^1.$ This implies $d^1$ is a split injection.

Therefore, $F_i \otimes_A F_i'$ is homotopy equivalent to a module $V$ which satisfies
\[ (P_i \otimes \; _iP) \oplus (P_i \otimes \; _iP) \oplus V \simeq (P_i \otimes \; _iP \otimes_A P_i \otimes \; _iP) \oplus A.\]
Since $dim_K(P_i \otimes_A \; _iP) = dim_K Hom_A(P_i,P_i) = 2,$ we obtain $V \simeq A$ as a left module and this finishes the proof. \qed

Denote by \[\mathcal{F}_i := F_i \otimes_A -.\] It follows from above that : \[ \mathcal{F}_i^{-1} = F_i' \otimes_A - .\]

\begin{Lemma} \label{imageFi} We have
\[ \mathcal{F}_i(P_j) =  \left\{ \begin{array}{cccccccc}
0 & \rightarrow & P_i & \rightarrow & 0 &               &   & \hbox{if }\; i = j \\
0 & \rightarrow & P_i & \rightarrow & P_j & \rightarrow & 0 &  \hbox{if }\; i \neq j \end{array} \right. \]
\[ \mathcal{F}_i^{-1}(P_j) = \left\{ \begin{array}{cccccccc}
  &             &  0  & \rightarrow & P_i & \rightarrow & 0 & \hbox{if }\; i = j \\
0 & \rightarrow & P_j & \rightarrow & P_i & \rightarrow & 0 & \hbox{if }\; i \neq j \end{array} \right. \] where $P_j$ is in degree 0.
\end{Lemma}

\noindent {\bf Proof} The images of $P_j$ under $\mathcal{F}_i$ and $\mathcal{F}_i^{-1}$ are \[\mathcal{F}_i(P_j) = ( 0 \rightarrow P_i \otimes \; _iP \otimes_A P_j \stackrel{\alpha_i \otimes id_{P_j}}{\longrightarrow} A \otimes_A P_j \rightarrow 0 ) \;\; \hbox{and}\]
 \[\mathcal{F}_i^{-1}(P_j) = ( 0 \rightarrow A \otimes_A P_j \stackrel{\beta_i \otimes id_{P_j}}{\longrightarrow} P_i \otimes \; _iP \otimes_A P_j \rightarrow 0 ).\]
When $i \neq j,$ $_iP \otimes_A P_j$ is generated by the path $(i \cdots j),$ hence
\[ \alpha_i \otimes id_{P_j}(x \otimes y \otimes_A z) = x(i \cdots j) \; \hbox{for all}\; x \otimes y \otimes_A z \in P_i \otimes \; _iP \otimes_A P_j, \]
\[ \beta_i \otimes id_{P_j}(1 \otimes_A x) = x(j \cdots i) \otimes (i \cdots j) \; \hbox{for all} \; x \in P_j.\]
Since $dim_K(_iP \otimes_A P_j) = dim_KHom_A(P_i,P_j) = 1$ we have $P_i \otimes \; _iP \otimes_A P_j \simeq P_i.$ The only map up to scalars between $P_i$ and $P_j$ is
\[ \mu_{ij} : P_i \rightarrow P_j \; \hbox{where} \; \mu_{ij}(x) = x(i \cdots j) \; \hbox{for all} \; x \in P_i.\]
Therefore, $\alpha_i \otimes id_{P_j} \simeq \mu_{ij},$ $\beta_i \otimes id_{P_j} \simeq \mu_{ji}$ and
\[\mathcal{F}_i(P_j) = (0 \rightarrow P_i \stackrel{\mu_{ij}}{\rightarrow} P_j \rightarrow 0) \; \hbox{and} \; \mathcal{F}_i^{-1}(P_j) = (0 \rightarrow P_j \stackrel{\mu_{ji}}{\rightarrow} P_i \rightarrow 0). \]
When $i = j,$
$ P_i \otimes \; _iP \otimes_A P_i \stackrel{\alpha_i \otimes id_{P_i}}{\longrightarrow} A \otimes_A P_i$ is surjective since for all $x \in P_i,$ \[(\alpha_i \otimes id_{P_i})(x \otimes (i) \otimes (i) ) = x.\]
Now if  $x \in P_i$, then one can check that \[\beta_i \otimes id_{P_i}(1 \otimes_A x) = x \otimes (i \cdots i) + x(i \cdots i) \otimes (i).\] Hence, if $x \neq 0,$ then $\beta_i \otimes id_{P_i}(1 \otimes_A x) \neq 0;$ therefore,  $A \otimes_A P_i \stackrel{\beta_i \otimes id_{P_i}}{\longrightarrow} P_i \otimes \; _iP \otimes_A P_i$ is injective.

 Hence, $\mathcal{F}_i(P_i)$ (resp. $\mathcal{F}_i^{-1}(P_i)$) has homology concentrated in degree 1 (resp. $-1$). As  $dim_K(_iP \otimes_A P_i) =  2,$ $\mathcal{F}_i(P_i) \simeq P_i[1]$ and $\mathcal{F}_i^{-1}(P_i) \simeq P_i[-1].$\qed

Denote by $\Delta$ the subgroup of $TrPic(A)$ generated by $F_1,
\cdots, F_n.$
\begin{Lemma} \label{enough}  For $F, F' \in \Delta,$ $F \simeq F'$ if and only if
\[ F \otimes_A P \simeq F' \otimes_A P \] for any  indecomposable projective module $P.$ \end{Lemma}

\noindent {\bf Proof} Let us assume that $F \otimes_A P \simeq F' \otimes_A P$ for any  indecomposable projective $A$-module $P.$
Denote by $B$ the Brauer tree algebra associated to a line with $n$ edges and multiplicity 1. Denote by $\Phi$ the isomorphism of groups \[ \Phi : TrPic(A) \rightarrow TrPic(B) \] and by $G$ and $G'$ the images of $F$ and $F'$ under $\Phi$ ($G := \Phi(F)$ and $G' := \Phi(F')).$
Denote by $X$ the associated two-sided tilting complex in $D^b(B \otimes A^{op})$ giving the derived equivalence between $A$ and $B,$ so we have \[G^{-1} \otimes_B X \simeq X \otimes_A F^{-1} \;\;\hbox{and}\;\; G' \otimes_B X \simeq X \otimes_A F'.\]
\[ \begin{array}{ccc} & \stackrel{F \otimes_A -}{\rightarrow} &  \\  D^b(A) & \stackrel{F^{-1} \otimes_A - }{\leftarrow} & D^b(A)  \\ \\ _{X \otimes_A -}  \downarrow \;\;\;\;\;\;\;\;\;\;\;\;\;\;& & \;\;\;\;\;\;\;\;\;\;\;\;\;\;\;\downarrow\;  _{ X \otimes_A -} \\ \\  D^b(B) & \stackrel{G \otimes_B -}{\rightarrow} & D^b(B)  \\  &  \stackrel{G^{-1} \otimes_B - }{\leftarrow} &  \end{array}\]
We will show that $G^{-1} \otimes_B Q \simeq G'^{-1} \otimes_B Q,$ or equivalently \[Q \simeq G' \otimes_B G^{-1} \otimes_B Q\]  for any  indecomposable projective $B$-module $Q.$ \\
Let $Q$ be an  indecomposable projective $B$-module and denote by $T$ the partial tilting complex of $A$ where $Q = \Phi(T) = X \otimes_A T$ and denote by $P$ the  indecomposable projective $A$-module where $P := F^{-1} \otimes_A T.$ Then we have
\begin{eqnarray*} G' \otimes_B G^{-1} \otimes_B Q & \simeq & G' \otimes_B G^{-1} \otimes_B X \otimes_A T \\
                                                  & \simeq & G' \otimes_B X \otimes_A F^{-1} \otimes_A T \\
                                                  & \simeq & G' \otimes_B X \otimes_A P \\
                                                  & \simeq & X \otimes_A F' \otimes_A P \\
                                                  & \simeq & X \otimes_A F \otimes_A P \;\; \hbox{by assumption} \\
                                                  & \simeq & X \otimes_A F \otimes_A F^{-1} \otimes_A T \\
                                                  & \simeq & X \otimes_A T \simeq Q. \end{eqnarray*}
Hence, by \cite[Remark 3]{RouZi}, using  the identity component $Out_0(B)$ of the outer automorphism group of $B$ we obtain $G^{-1} \simeq G'^{-1},$ which implies $G \simeq G'.$ And since $Out_0$ is invariant under derived equivalence (see \cite{Huisgen-Saorin}) we obtain  $F \simeq F'.$
 \qed

\section{Braid action}

Denote by $B(K_n)$ the braid group whose associated Coxeter group has Coxeter graph given by the complete graph on $n$ vertices, generated by $a_1, \cdots, a_n$ with relations \[a_i a_j a_i = a_j a_i a_j\] for all $1 \le i,j \le n, i \neq j.$ We will define an action of $B(K_n)$ on $D^b(A).$

\begin{Theorem} \label{action} There is a group homomorphism \[ \varphi : B(K_n) \rightarrow TrPic(A) \] where $\varphi(a_i) = F_i.$
\end{Theorem}

\noindent {\bf Proof} By Lemma \ref{enough}, it is enough to show that $\mathcal{F}_i \mathcal{F}_j \mathcal{F}_i (P) \simeq \mathcal{F}_j \mathcal{F}_i \mathcal{F}_j(P)$ for every indecomposable projective module $P.$ We will show that for every indecomposable projective module $P$
\[ \mathcal{F}_i^{-1} \mathcal{F}_j \mathcal{F}_i (P) \simeq \mathcal{F}_j \mathcal{F}_i \mathcal{F}_j^{-1} (P). \]  \\

{\bf Claim}
{\footnotesize
\[ \left( \begin{array}{c} P_1 \\ P_2 \\ \vdots \\ P_{i-1} \\ P_i \\ P_{i+1} \\ \vdots \\ P_{j-1} \\ P_j \\ P_{j+1} \\ \vdots \\ P_n \end{array} \right) \stackrel{\mathcal{F}_i}{\longrightarrow}
\left( \begin{array}{ccc} P_i & \rightarrow & P_1 \\ P_i & \rightarrow & P_2 \\& \vdots & \\ P_i & \rightarrow & P_{i-1} \\ P_i & \rightarrow  & 0 \\ P_i & \rightarrow & P_{i+1} \\& \vdots & \\ P_i & \rightarrow & P_{j-1} \\ P_i & \rightarrow & P_j \\ P_i & \rightarrow & P_{j+1} \\& \vdots & \\ P_i & \rightarrow & P_n \end{array} \right)
\stackrel{\mathcal{F}_j}{\longrightarrow}
\left( \begin{array}{ccccc} & & P_i & \rightarrow & P_1 \\ & &  P_i & \rightarrow & P_2 \\& & & \vdots  & \\& & P_i & \rightarrow & P_{i-1} \\P_j & \rightarrow & P_i & \rightarrow  & 0 \\& &  P_i & \rightarrow & P_{i+1} \\& & & \vdots  & \\ & &  P_i & \rightarrow & P_{j-1} \\ & & P_i & \rightarrow & 0\\ & & P_i & \rightarrow & P_{j+1} \\& &  & \vdots  & \\& & P_i & \rightarrow & P_n \end{array} \right)
\stackrel{\mathcal{F}_i^{-1}}{\longrightarrow}
\left( \begin{array}{cc}  & P_1 \\   & P_2 \\ &  \vdots \\  & P_{i-1} \\P_j  \rightarrow P_i \rightarrow  & P_i \\ & P_{i+1} \\ &  \vdots \\ &  P_{j-1} \\ & P_i\\ & P_{j+1} \\ & \vdots \\&  P_n  \end{array} \right)\]} and
{\footnotesize
\[ \left( \begin{array}{c} P_1 \\ P_2 \\ \vdots \\ P_{i-1} \\ P_i \\ P_{i+1} \\ \vdots \\ P_{j-1} \\ P_j \\ P_{j+1} \\ \vdots \\ P_n \end{array} \right) \stackrel{\mathcal{F}_j^{-1}}{\longrightarrow}
\left( \begin{array}{ccc} P_1 & \rightarrow & P_j \\ P_2 & \rightarrow & P_j \\& \vdots & \\ P_{i-1} & \rightarrow & P_j \\ P_i & \rightarrow  & P_j \\ P_{i+1} & \rightarrow & P_j \\& \vdots & \\ P_{j-1} & \rightarrow & P_j \\  & & P_j \\ P_{j+1} & \rightarrow & P_j \\& \vdots & \\ P_n & \rightarrow & P_j \end{array} \right)
\stackrel{\mathcal{F}_i}{\longrightarrow}
\left( \begin{array}{ccccc} & & P_1 & \rightarrow & P_j \\ & &  P_2 & \rightarrow & P_j \\& & & \vdots  & \\& & P_{i-1} & \rightarrow & P_j \\P_i & \rightarrow & P_i & \rightarrow  & P_j \\& &  P_{i+1} & \rightarrow & P_j \\& & & \vdots  & \\ & &  P_{j-1} & \rightarrow & P_j \\ & & P_i & \rightarrow & P_j\\ & & P_{j+1} & \rightarrow & P_j \\& &  & \vdots  & \\& & P_n & \rightarrow & P_j \end{array} \right)
\stackrel{\mathcal{F}_j}{\longrightarrow}
\left( \begin{array}{cc}  & P_1 \\  & P_2 \\ &  \vdots \\  & P_{i-1} \\P_j  \rightarrow P_i \rightarrow  & P_i \\ & P_{i+1} \\ &  \vdots \\  &   P_{j-1} \\  & P_i\\ &   P_{j+1} \\ & \vdots \\& P_n  \end{array} \right)\]}
\begin{equation} \label{FiFjFi}
\end{equation}

\noindent {\bf Proof of Claim :}\\
For all $1 \le i \le n$, denote by $\delta_i : P_i \rightarrow P_i$ where $\delta_i((i)) = (i),$  and for all $1 \le i,j \le n$ by $\mu_{ij} : P_i \rightarrow P_j,$ the only map up to a scalars between $P_i$ and $P_j$ where
\[\mu_{ij}(x) = x(i \cdots j) \; \hbox{for all} \; x \in P_i.\] \\
Since $\mathcal{F}_j$ preserves cones, $\mathcal{F}_j(P_i \rightarrow P_j)$   is isomorphic to
{\small \begin{eqnarray*}  Cone (\mathcal{F}_j(P_i)  \rightarrow \mathcal{F}_j(P_j)) & = &
Cone\left(\begin{array}{cccc} \mathcal{F}_j(P_i) = & P_j \otimes \;_jP \otimes_A P_i & \rightarrow & A \otimes_A P_i   \\
                                                &  x \otimes (j \cdots i)          &             & x(j \cdots i) \\ \\
                             \downarrow         & \downarrow                      &             & \downarrow     \\ \\
                           \mathcal{F}_j(P_j) =  &   P_j \otimes \;_jP \otimes_A P_j & \rightarrow & A \otimes_A P_j   \\
                                                &  x \otimes (j \cdots j)          &             & x(j \cdots j)   \end{array}\right) \\ &  \simeq & Cone\left(\begin{array}{ccc} P_j & \stackrel{\mu_{ji}}{\rightarrow} & P_i \\
                                            x  &                                  & x(j \cdots i) \\
                                          \;\;\; \downarrow  \;_{\delta_j}                                           \\
                                           P_j &                                  &               \\
                                           x   &                                  &           \end{array} \right) \\ & \simeq & P_i[1].\end{eqnarray*}}
For $k \neq j,$
 $\mathcal{F}_j(P_i \rightarrow P_k)$ is isomorphic to
{\small \[
Cone\left(\begin{array}{cccc} \mathcal{F}_j(P_i) = & P_j \otimes \;_jP \otimes_A P_i & \rightarrow & A \otimes_A P_i   \\
                                                &  x \otimes (j \cdots i)          &             & x(j \cdots i) \\ \\
                             \downarrow         & \downarrow                      &             & \downarrow \;_{\mu_{ik}}    \\ \\
                           \mathcal{F}_j(P_k) =  &   P_j \otimes \;_jP \otimes_A P_k & \rightarrow & A \otimes_A P_k   \\
                                                &  x \otimes (j \cdots k)          &             & x(j \cdots k)   \end{array}\right)   \simeq  Cone\left(\begin{array}{ccc} P_j & \stackrel{\mu_{ji}}{\rightarrow} & P_i \\
                                            x  &                                  & x(j \cdots i) \\
                                         \;\;\;  \downarrow \; _{\delta_j}& & \;\;\; \downarrow \; _{\mu_{ik}}                                            \\
                                           P_j &  \stackrel{\mu_{jk}}{\rightarrow}   & P_k              \\
                                           x   &                                  & x(j \cdots k)          \end{array} \right) \]
\begin{eqnarray*} & \simeq &  \begin{array}{ccccc} P_j & \stackrel{\mu_{ji}}{\rightarrow} & P_i & & \\ & _{\delta_j} \searrow \;\;\;\;\;\; & \oplus &\;\;\;\;\;\;\searrow \;_{\mu_{ik}} & \\ & & P_j & \stackrel{\mu_{jk}}{\rightarrow} & P_k \end{array} \\ & \simeq & Cone( Cone (P_j \rightarrow P_i \oplus P_j) \rightarrow P_k ) \\ & \simeq &   Cone \left( \begin{array}{c} Cone \left( \begin{array}{ccc} (P_j & \stackrel{\mu_{ji}}{\rightarrow} & P_i) \\ _{\delta_j} \downarrow & &  \\ P_j & & \end{array} \right) [-1] \;\;\; \stackrel{(\mu_{ik},\mu_{jk})}{\longrightarrow} P_k \end{array} \right) \\& \simeq &  P_i \stackrel{\mu_{ik}}{\longrightarrow} P_k.
\end{eqnarray*}}
The complex $\mathcal{F}_i^{-1}(P_i \rightarrow P_k)$ ($k \neq i$) is isomorphic to
{\small \[
Cone\left(\begin{array}{cccc} \mathcal{F}^{-1}_i(P_i) = & A \otimes_A P_i  & \rightarrow & P_i \otimes \;_iP \otimes_A P_i  \\                                                     &  x               &             & x(i \cdots i) \otimes (i) + x \otimes (i \cdots i) \\ \\
                             \downarrow         & \downarrow                      &             & \downarrow    \\ \\
                           \mathcal{F}^{-1}_i(P_k) =  & A \otimes_A P_k  & \rightarrow &  P_i \otimes \;_iP \otimes_A P_k   \\
                                                &  x(i \cdots k)          &             & x(i \cdots i) \otimes (i \cdots k)  \end{array}\right)\] \[ \simeq
Cone \left( \begin{array}{ccc} & & P_i  \\ & & \downarrow \; _{\delta_i}\\ (P_k & \stackrel{\mu_{ki}}{\rightarrow} & P_i) \end{array} \right) \simeq P_k \;\; (P_i \;\hbox{lies in degree}\; -1).\]}
The complex $\mathcal{F}_i^{-1}(P_j \rightarrow P_i \rightarrow 0)$  is isomorphic to
{\small \[
Cone\left(\begin{array}{cccc} \mathcal{F}^{-1}_i(P_j) = & A \otimes_A P_j  & \rightarrow & P_i \otimes \;_iP \otimes_A P_j  \\                                                     &  x               &             & x(j \cdots i) \otimes (i \cdots j) \\ \\
                             \downarrow         & \downarrow                      &             & \downarrow    \\ \\
                           \mathcal{F}^{-1}_i(P_i) =  & A \otimes_A P_i  & \rightarrow &  P_i \otimes \;_iP \otimes_A P_i   \\
                                                &  x(j \cdots i)          &             & x(j \cdots i) \otimes (i \cdots i)  \end{array}\right)\]
\[ \simeq Cone \left( \begin{array}{ccc} (P_j & \stackrel{\mu_{ji}}{\rightarrow} & P_i) \\ x & & x(j \cdots i) \\ & & \;\;\;\;\downarrow \;_{\delta_{i}} \\ & & P_i \\ & & x(j \cdots i) \end{array} \right) \simeq P_j \stackrel{\mu_{ji}}{\longrightarrow} P_i \stackrel{\delta_i}{\longrightarrow} P_i,\]}  and the complex $\mathcal{F}_i(P_i \rightarrow P_j)$ is isomorphic to
{\small  \[
Cone\left(\begin{array}{cccc} \mathcal{F}_i(P_i) = & P_i \otimes \;_iP \otimes_A P_i & \rightarrow & A \otimes_A P_i   \\
                                                &   x \otimes (i) + x \otimes (i \cdots i)         &             & x(i) \\ \\
                             \downarrow         & \downarrow                      &             & \downarrow \; _{\mu_{ij}}    \\ \\
                           \mathcal{F}_i(P_j) =  &   P_i \otimes \;_iP \otimes_A P_j & \rightarrow & A \otimes_A P_j  \\
                                                &  x \otimes (i \cdots j)          &             & x(i \cdots j)   \end{array}\right)\] \[\simeq Cone \left( \begin{array}{ccc} P_i & & \\ x & & \\ _{\delta_i} \downarrow \;\;\;\; & & \\ (P_i & \stackrel{\mu_{ij}}{\rightarrow} & P_j) \\ x & & x(i \cdots j) \end{array} \right) \simeq P_i \stackrel{\delta_i}{\longrightarrow} P_i \stackrel{\mu_{ij}}{\longrightarrow} P_j\]}
And finally, in a similar way, the complex $\mathcal{F}_j(P_i \rightarrow P_i \rightarrow P_j)$ (where $P_j$ is in degree $-1$) is isomorphic to
{\small \begin{eqnarray*}\mathcal{F}_j(Cone (P_i \rightarrow (Cone (P_i \rightarrow P_j))[-1])) & \simeq & Cone (\mathcal{F}_j(P_i) \rightarrow (\mathcal{F}_j(P_i \rightarrow P_j))[-1]) \\ & \simeq & Cone \left( \begin{array}{ccc}( P_j & \stackrel{\mu_{ji}}{\rightarrow} & P_i) \\ & & _{\delta_i} \downarrow \;\;\;\; \\  & & P_i \end{array} \right) \\ & \simeq& P_j \stackrel{\mu_{ij}}{\longrightarrow} P_i \stackrel{\delta_i}{\longrightarrow} P_i. \end{eqnarray*}}
This finishes the proof. \qed

\section{Categorification of a braid group homomorphism}

Let $B$ be a Brauer tree algebra associated to a line with $n$ edges and no exceptional vertex.
Denote by $B(A_n)$ the Artin braid group on $n$ strings generated by $\sigma_1,\cdots,\sigma_n$ satisfying the relations
\[ \sigma_i \sigma_{i+1} \sigma_i = \sigma_{i+1} \sigma_i \sigma_{i+1} \; \hbox{and} \; \sigma_i \sigma_j = \sigma_j \sigma_i \; \hbox{if} \; |i - j| > 1. \]
In \cite{RouZi} Rouquier and Zimmermann obtained an action of $B(A_n)$ on  derived category of $B.$

\begin{Theorem} (Rouquier and Zimmermann \cite{RouZi})\\ There is a group homomorphism \begin{eqnarray*} \psi_n :&  B(A_n) & \rightarrow  TrPic(B) \\
                           & \sigma_i & \mapsto  R_i \end{eqnarray*}
where $R_i$ is the two-sided tilting complex $ 0 \rightarrow Q_i \otimes Hom_K(Q_i,K) \rightarrow B \rightarrow 0$ and $Q_1, \cdots Q_n$ are the indecomposable projective modules of $B.$ Moreover $\psi_2$ is injective.
\end{Theorem}

In \cite{Khovanov-Seidel} Khovanov and Seidel independently discovered this homomorphism $\psi$ for similar algebras in a very different context.

\begin{Theorem} (Khovanov and Seidel \cite{Khovanov-Seidel})\\ The map $\psi_n$ is injective for all $n.$
\end{Theorem}

For simplicity, let us denote by $\psi := \psi_n.$
Denote by $\mathcal{R}_i$ the corresponding functor $\mathcal{R}_i := R_i \otimes_B -.$ By \cite[Lemma 4.2]{RouZi}, the images of the indecomposable projective modules $Q_j$ of $B$ under $\mathcal{R}_i$ are
{\small \[ \mathcal{R}_i(Q_j) =  \left\{ \begin{array}{cccccccc}
0 & \rightarrow & Q_i & \rightarrow & 0 &               &   & \hbox{if }\; i = j \\
0 & \rightarrow & Q_i & \rightarrow & Q_j & \rightarrow & 0 &  \hbox{if }\; |i - j| = 1 \\
  &             &  0  & \rightarrow & Q_j & \rightarrow & 0 & \hbox {if} \; |i - j| > 1   \end{array} \right. \]}
 where $Q_j$ is in degree 0. \\

In the next subsections we will provide a link between the $B(K_n)$, $B(A_n)$, the algebras $A$ and $B$ which will show that the action defined in Theorem \ref{action} is unfaithful.

\subsection{Braid group homomorphism}

Define
\[ c_n  :=  \sigma_n \;\; \hbox{and} \;\;  c_k  :=   \sigma_k^{-1} c_{k+1} \sigma_k \;\; \hbox{for all} \;\; 1 \le k \le n-1. \] Define a mapping
\begin{eqnarray*} \eta : & B(K_n) & \rightarrow B(A_n) \\
                         & a_i    & \mapsto c_i  \end{eqnarray*}

\begin{Proposition} \label{eta} The mapping $\eta$ is a surjective group homomorphism. \end{Proposition}

First we will need the following Lemma :

\begin{Lemma} \label{need} $c_j \sigma_i^{-1} = \sigma_i^{-1}c_j$  for all $1 \le j \le n-1,\; 1 \le i \le n,\; i \neq j-1,\; i \neq j$ \end{Lemma}

\noindent {\bf Proof} If $i > j,$ we have
{\small \begin{eqnarray*} c_j \sigma_i^{-1} & = & \sigma_j^{-1} \cdots \sigma_{n-1}^{-1} \sigma_n \sigma_{n-1} \cdots \sigma_j \sigma_i^{-1} \\
                               & = & \sigma_j^{-1} \cdots \sigma_{n-1}^{-1} \sigma_n \sigma_{n-1} \cdots \sigma_{i+1} \sigma_i \sigma_{i-1} \sigma_i^{-1} \sigma_{i-2} \cdots \sigma_j  \\             & = &  \sigma_j^{-1} \cdots \sigma_{n-1}^{-1} \sigma_n \sigma_{n-1} \cdots \sigma_{i+1} \sigma_{i-1}^{-1} \sigma_i \sigma_{i-1} \sigma_{i-2} \cdots \sigma_j \\          & = &  \sigma_j^{-1} \cdots \sigma_{i-1}^{-1} \sigma_i^{-1} \sigma_{i-1}^{-1} \sigma_{i+1}^{-1} \cdots \sigma_{n-1}^{-1} \sigma_n \sigma_{n-1} \cdots  \sigma_{i+1} \sigma_i \sigma_{i-1} \sigma_{i-2} \cdots \sigma_j \\
 & = &  \sigma_j^{-1} \cdots \sigma_i^{-1} \sigma_{i-1}^{-1} \sigma_i^{-1} \sigma_{i+1}^{-1} \cdots \sigma_{n-1}^{-1} \sigma_n \sigma_{n-1} \cdots \sigma_j \\
 & = & \sigma_i^{-1} c_j \end{eqnarray*}}
If $i < j,$ $i \neq j-1,$
{\small \begin{eqnarray*} c_j \sigma_i^{-1} &= & \sigma_j^{-1} \cdots \sigma_{n-1}^{-1} \sigma_n \sigma_{n-1} \cdots \sigma_j \sigma_i^{-1} \\
                               & = & \sigma_i^{-1}  \sigma_j^{-1} \cdots \sigma_{n-1}^{-1} \sigma_n \sigma_{n-1} \cdots \sigma_j \\
                               & = & \sigma_i^{-1} c_j \;\;\;\;  \end{eqnarray*}} \qed

\noindent {\bf Proof of Proposition \ref{eta}} We need to show \begin{enumerate} \item $\sigma_n c_k \sigma_n = c_k \sigma_n c_k$ for all $1 \le k \le n-1.$
\item $c_i c_j c_i = c_j c_i c_j$ for all $1 \le i,j \le n-1,$ $i \neq j.$
\item $\sigma_i$ is in the image of $\eta$ for all $1 \le i \le n.$ \end{enumerate}

\begin{enumerate}
\item We will show this using induction. For $k = n-1,$
{\small \begin{eqnarray*} \sigma_n c_{n-1} \sigma_n & = & \sigma_n \sigma_{n-1}^{-1} \sigma_n \sigma_{n-1} \sigma_n \\
                                  & = & \sigma_n \sigma_n \sigma_{n-1} \sigma_n^{-1} \sigma_n \\
                                  & = & \sigma_n \sigma_n \sigma_{n-1}. \end{eqnarray*}
\begin{eqnarray*} c_{n-1} \sigma_n c_{n-1} & = & \sigma_{n-1}^{-1} \sigma_n \sigma_{n-1} \sigma_n \sigma_{n-1}^{-1} \sigma_n \sigma_{n-1} \\                 & = & \sigma_{n-1}^{-1} \sigma_n \sigma_n^{-1} \sigma_{n-1} \sigma_n \sigma_n \sigma_{n-1} \\                     & = & \sigma_n \sigma_n \sigma_{n-1}. \end{eqnarray*}}
Assume that $\sigma_n c_j \sigma_n = c_j \sigma_n c_j;$ we will show that $\sigma_n c_{j-1} \sigma_n = c_{j-1} \sigma_n c_{j-1}.$
{\small \begin{eqnarray*} \sigma_n c_{j-1} \sigma_n & = & \sigma_n \sigma_{j-1}^{-1} c_j \sigma_{j-1} \sigma_n \\
                                  & = & \sigma_{j-1}^{-1} \sigma_n c_j \sigma_n \sigma_{j-1} \\
                                  & = & \sigma_{j-1}^{-1} c_j \sigma_n c_j \sigma_{j-1} \\
                                  & = & \sigma_{j-1}^{-1} c_j \sigma_{j-1} \sigma_n \sigma_{j-1}^{-1} c_j \sigma_{j-1} \\                 & = & c_{j-1} \sigma_n c_{j-1}. \end{eqnarray*}}
 \item For $i = n-1,$  if $j = n-2, $
{\small \begin{eqnarray*}
c_{n-1} c_{n-2} c_{n-1} & = &  \sigma_{n-1}^{-1} \sigma_n \sigma_{n-1} \sigma_{n-2}^{-1} \sigma_{n-1}^{-1} \sigma_n \sigma_{n-1} \sigma_{n-2} \sigma_{n-1}^{-1} \sigma_n \sigma_{n-1} \\
& = &  \sigma_n \sigma_{n-1} \sigma_n^{-1} \sigma_{n-2}^{-1} \sigma_n \sigma_{n-1} \sigma_n^{-1} \sigma_{n-2} \sigma_n \sigma_{n-1} \sigma_n^{-1} \\
& = & \sigma_n \sigma_{n-1} \sigma_{n-2}^{-1} \sigma_{n-1} \sigma_{n-2} \sigma_{n-1} \sigma_n^{-1}\\
& = & \sigma_n \sigma_{n-1} \sigma_{n-2}^{-1} \sigma_{n-2} \sigma_{n-1} \sigma_{n-2} \sigma_n^{-1} \\
& = & \sigma_n \sigma_{n-1} \sigma_{n-1} \sigma_{n-2} \sigma_n^{-1}. \end{eqnarray*}
\begin{eqnarray*}
c_{n-2} c_{n-1} c_{n-2} & = &  \sigma_{n-2}^{-1} \sigma_{n-1}^{-1} \sigma_n \sigma_{n-1} \sigma_{n-2} \sigma_{n-1}^{-1} \sigma_n \sigma_{n-1} \sigma_{n-2}^{-1} \sigma_{n-1}^{-1} \sigma_n \sigma_{n-1} \sigma_{n-2} \\
& = &  \sigma_{n-2}^{-1} \sigma_n \sigma_{n-1} \sigma_n^{-1} \sigma_{n-2}  \sigma_n \sigma_{n-1} \sigma_n^{-1} \sigma_{n-2}^{-1} \sigma_n \sigma_{n-1} \sigma_n^{-1} \sigma_{n-2} \\
& = & \sigma_{n-2}^{-1} \sigma_n \sigma_{n-1} \sigma_{n-2} \sigma_{n-1} \sigma_{n-2}^{-1} \sigma_{n-1} \sigma_n^{-1} \sigma_{n-2}  \\
& = & \sigma_{n-2}^{-1} \sigma_n \sigma_{n-2} \sigma_{n-1} \sigma_{n-2} \sigma_{n-2}^{-1} \sigma_{n-1} \sigma_{n-2} \sigma_n^{-1} \\
& = & \sigma_n \sigma_{n-1} \sigma_{n-1} \sigma_{n-2} \sigma_n^{-1}. \end{eqnarray*}}
If $j \neq n-2,$ $c_{n-1} c_j c_{n-1} =$
{\small \begin{eqnarray*}
 & = & \sigma_{n-1}^{-1} \sigma_n \sigma_{n-1} (\sigma_j^{-1} \sigma_{j+1}^{-1} \cdots \sigma_{n-3}^{-1} \sigma_{n-2}^{-1} \sigma_{n-1}^{-1} \sigma_n \sigma_{n-1} \sigma_{n-2} \sigma_{n-3} \cdots \sigma_{j+1} \sigma_j) \sigma_{n-1}^{-1} \sigma_n \sigma_{n-1} \\
 & = & (\sigma_j^{-1} \sigma_{j+1}^{-1} \cdots \sigma_{n-3}^{-1}) \sigma_{n-1}^{-1} \sigma_n \sigma_{n-1} (\sigma_{n-2}^{-1} \sigma_{n-1}^{-1} \sigma_n \sigma_{n-1} \sigma_{n-2}) \sigma_{n-1}^{-1} \sigma_n \sigma_{n-1} (\sigma_{n-3} \cdots \sigma_j) \\
 & = &  \sigma_j^{-1} \sigma_{j+1}^{-1} \cdots \sigma_{n-3}^{-1} c_{n-1} c_{n-2} c_{n-1} \sigma_{n-3} \cdots \sigma_j \\
 & = &  \sigma_j^{-1} \sigma_{j+1}^{-1} \cdots \sigma_{n-3}^{-1} c_{n-2} c_{n-1} c_{n-2} \sigma_{n-3} \cdots \sigma_j, \end{eqnarray*}
$c_j c_{n-1} c_j =$
\begin{eqnarray*}
 & = & \sigma_j^{-1} \sigma_{j+1}^{-1} \cdots \sigma_{n-1}^{-1} \sigma_n \sigma_{n-1} \cdots \sigma_j   \sigma_{n-1}^{-1} \sigma_n \sigma_{n-1}  \sigma_j^{-1} \sigma_{j+1}^{-1} \cdots \sigma_{n-1}^{-1} \sigma_n \sigma_{n-1} \cdots \sigma_j \\
 & = &  \sigma_j^{-1} \sigma_{j+1}^{-1} \cdots \sigma_{n-3}^{-1}  \sigma_{n-2}^{-1} \sigma_{n-1}^{-1} \sigma_n \sigma_{n-1} \sigma_{n-2} \sigma_{n-1}^{-1} \sigma_n \sigma_{n-1} \sigma_{n-2}^{-1} \sigma_{n-1}^{-1} \sigma_n \sigma_{n-1} \sigma_{n-2} \sigma_{n-3} \cdots \sigma_j \\
 & = &  \sigma_j^{-1} \sigma_{j+1}^{-1} \cdots \sigma_{n-3}^{-1} c_{n-2} c_{n-1} c_{n-2} \sigma_{n-3} \cdots \sigma_j. \end{eqnarray*}}

Now assume that $c_j c_k c_j = c_k c_j c_k$ for all $k > i.$ We will show that  \[c_j c_i c_j = c_i c_j c_i.\]
If $j \neq i+1$
{\small \begin{eqnarray*} c_j c_i c_j & = & c_j \sigma_i^{-1} c_{i+1} \sigma_i c_j \\
                              & = & \sigma_i^{-1} c_j c_{i+1} c_j \sigma_i \;\; \hbox{by Lemma \ref{need}}\\
                              & = & \sigma_i^{-1} c_{i+1} c_j c_{i+1} \sigma_i \\
                              & = & \sigma_i^{-1} c_{i+1} \sigma_i c_j \sigma_i^{-1} c_{i+1} \sigma_i \\
                              & = & c_i c_j c_i. \end{eqnarray*}}
If $j = i+1$
{\small \begin{eqnarray*} c_{i+1} c_i c_{i+1} & = & c_{i+1} \sigma_i^{-1} c_{i+1} \sigma_i c_{i+1} \\
                                      & = & c_{i+1} \sigma_i^{-1} \sigma_{i+1}^{-1} c_{i+2} \sigma_{i+1} \sigma_i c_{i+1} \\
                                      & = & \sigma_{i+1}^{-1} c_{i+2} \sigma_{i+1}  \sigma_i^{-1} \sigma_{i+1}^{-1} c_{i+2} \sigma_{i+1} \sigma_i \sigma_{i+1}^{-1} c_{i+2} \sigma_{i+1} \\          & = & \sigma_{i+1}^{-1} c_{i+2} \sigma_i^{-1} \sigma_{i+1}^{-1} \sigma_i c_{i+2} \sigma_i^{-1} \sigma_{i+1} \sigma_i  c_{i+2} \sigma_{i+1} \\                  & = & \sigma_{i+1}^{-1} \sigma_i^{-1} c_{i+2} \sigma_{i+1}^{-1} c_{i+2} \sigma_{i+1} \sigma_i c_{i+2} \sigma_{i+1} \\                                & = & \sigma_{i+1}^{-1} \sigma_i^{-1} c_{i+2} c_{i+1} c_{i+2} \sigma_i \sigma_{i+1} \\
                                      & = & \sigma_{i+1}^{-1} \sigma_i^{-1} c_{i+1} c_{i+2} c_{i+1} \sigma_i \sigma_{i+1} \\
                                      & = & \sigma_{i+1}^{-1} \sigma_i^{-1} \sigma_{i+1}^{-1} c_{i+2} \sigma_{i+1}  c_{i+2} \sigma_{i+1}^{-1} c_{i+2} \sigma_{i+1} \sigma_i \sigma_{i+1}\\            & = & \sigma_i^{-1} \sigma_{i+1}^{-1} \sigma_i^{-1} c_{i+2} \sigma_{i+1} c_{i+2} \sigma_{i+1}^{-1} c_{i+2} \sigma_i \sigma_{i+1} \sigma_i \\                  & = & \sigma_i^{-1} \sigma_{i+1}^{-1} c_{i+2} \sigma_i^{-1} \sigma_{i+1} c_{i+2} \sigma_{i+1}^{-1} \sigma_i c_{i+2} \sigma_{i+1} \sigma_i \\                   & = & \sigma_i^{-1} \sigma_{i+1}^{-1} c_{i+2} \sigma_i^{-1} \sigma_{i+1} \sigma_i c_{i+2} \sigma_i^{-1} \sigma_{i+1}^{-1} \sigma_i c_{i+2} \sigma_{i+1} \sigma_i \\      & = & \sigma_i^{-1} \sigma_{i+1}^{-1} c_{i+2} \sigma_{i+1} \sigma_i \sigma_{i+1}^{-1} c_{i+2} \sigma_{i+1} \sigma_i^{-1} \sigma_{i+1}^{-1} c_{i+2} \sigma_{i+1} \sigma_i \\
                                     & = & c_i c_{i+1} c_i. \end{eqnarray*}}

\item We will show that \[c_k c_{k+1} c_k^{-1} = c_{k+1}^{-1} c_k c_{k+1} = \sigma_k\] for all $1 \le k \le n-1.$
For $k = n-1,$ {\small \begin{eqnarray*} c_{n-1} c_n c_{n-1}^{-1} & = & \sigma_{n-1}^{-1} \sigma_n \sigma_{n-1} \sigma_n  \sigma_{n-1}^{-1} \sigma_n^{-1} \sigma_{n-1} \\                                                   & = & \sigma_n \sigma_{n-1} \sigma_n^{-1} \sigma_n \sigma_{n-1}^{-1} \sigma_n^{-1} \sigma_{n-1} \\
                                                          & = & \sigma_{n-1}. \end{eqnarray*}}
Assume $c_{k+1} c_{k+2} = c_{k+2} \sigma_{k+1}.$ We will show that $c_k c_{k+1} = c_{k+1} \sigma_k.$
{\small \begin{eqnarray*} c_k c_{k+1} & = & \sigma_k^{-1} c_{k+1} \sigma_k c_{k+1} \\
                              & = & \sigma_k^{-1} c_{k+2} \sigma_{k+1} c_{k+2}^{-1} \sigma_k c_{k+2} \sigma_{k+1} c_{k+2}^{-1} \\
                              & = &  \sigma_k^{-1} c_{k+2} \sigma_{k+1} \sigma_k  \sigma_{k+1} c_{k+2}^{-1} \\
                              & = &  \sigma_k^{-1} c_{k+2} \sigma_k  \sigma_{k+1}  \sigma_k c_{k+2}^{-1} \\
                              & = &   c_{k+2} \sigma_{k+1}  \sigma_k c_{k+2}^{-1} \\
                              & = & c_{k+1} c_{k+2} \sigma_k c_{k+2}^{-1} \\
                              & = & c_{k+1} \sigma_k. \;\;\;\;\;\; \end{eqnarray*} } \qed
\end{enumerate}

\begin{Proposition} \label{Rem} The homomorphism $\eta$ is not injective. \end{Proposition}

\noindent {\bf Proof} We have $\sigma_1 \sigma_n = \sigma_n \sigma_1$
but $a_1 a_2 a_1^{-1} a_n \neq a_n a_1 a_2 a_1^{-1}$
as we will show using the geometric representation of the Coxeter group of $B(K_n).$ \\

Denote by $W(K_n)$ the Coxeter group of $B(K_n).$ Denote by $d_i$ the image of $a_i$ in $W(K_n).$
By \cite[section 5.3]{Humphreys}, $W(K_n)$ acts on an $n$-dimensional vector space $V$ with basis $\{v_1, \cdots, v_n\}.$\\

The action of $W(K_n)$ is defined by $d_i \mapsto f_i : V \rightarrow V$ where
 \[ f_i(v_i) =  -v_i \; \hbox{and} \]
\[ f_i(v_j) =  v_i + v_j.\]
Now if $a_1 a_2 a_1^{-1} a_n a_1 a_2^{-1} a_1^{-1} a_n^{-1} = e,$ then \[d_1 d_2 d_1 d_n d_1 d_2 d_1 d_n = e \;\; \hbox{and} \;\; f_1 f_2 f_1 f_n f_1 f_2 f_1 f_n = id_V;\] but
 \begin{eqnarray*} f_1 f_2 f_1 f_n f_1 f_2 f_1 f_n(v_n) & = &  f_1 f_2 f_1 f_n f_1 f_2 f_1(-v_n)\\
                                                        & = &  f_1 f_2 f_1 f_n f_1 f_2 (-v_1 -v_n)\\
                                                        & = &  f_1 f_2 f_1 f_n f_1 (-2v_2 -v_1 -v_n) \\
                                                        & = &  f_1 f_2 f_1 f_n (-2v_1 -2 v_2 -v_n) \\
                                                        & = &  f_1 f_2 f_1 (-3v_n -2 v_1 -2 v_2) \\
                                                        & = &  f_1 f_2 (-3v_1 - 3v_n - 2v_2) \\
                                                        & = &  f_1 (-4v_2 - 3v_1 - 3v_n) \\
                                                        & = &  (-4v_1 - 4v_2 - 3v_n) \neq v_n, \;\; \hbox{a contradiction}. \end{eqnarray*}    \qed

\subsection{Categorification}

Let $T$ be the direct sum of the following complexes of projective $B$-modules.

\[ \begin{array}{ccccccccccccccc}
T_1 & : & 0 & \rightarrow & Q_n & \rightarrow & Q_{n-1} & \rightarrow & \cdots & \rightarrow & Q_2 & \rightarrow & Q_1 & \rightarrow & 0 \\
T_2 & : & 0 & \rightarrow & Q_n & \rightarrow & Q_{n-1} & \rightarrow & \cdots & \rightarrow & Q_2 & \rightarrow & 0 &  &  \\
\vdots    &   &   &             &     &             & \vdots  &             &        &             &     &             &   &  &  \\
T_{n-1} & : & 0 & \rightarrow & Q_n & \rightarrow & Q_{n-1} & \rightarrow & 0 & & & & & & \\
T_n     & : & 0 & \rightarrow & Q_n & \rightarrow & 0 & & & & & & & & \end{array} \]
where $Q_i$ is in degree $i-1.$ \\

By \cite[Theorem 4.2]{Rickarddercatstabeq}, using the stable equivalence between $End_{D^b(B)}(T)$ and $A$ and the result of Gabriel and Riedtmann \cite{Gabriel-Riedtmann},  $T$ is a tilting complex giving a derived equivalence between $B$ and $A.$ (We note that by computing the endomorphism ring of $T$  directly with the method in \cite{Muchtadi-Alamsyah}, we can show that $End_{D^b(B)}(T) \simeq A$). \\

Denote by $G$ the unique two-sided tilting complex of $(A \otimes B^{op})$-modules associated to $T$ (such a complex is unique up to isomorphism in $D^b(A \otimes B^{op})$ by a theorem of Keller \cite{Kellerremark}).  Denote by $\mathcal{G}$ the corresponding functor \[ \mathcal{G} = G \otimes_B -. \]
It is easy to see that \begin{equation} \label{GTi} \mathcal{G}(T_i) = P_i \;\; \hbox{for all} \;\;1 \le i \le n. \end{equation}

We will compute the image of $R_i$ under this functor.
\[ \begin{array}{ccc} D^b(B) & \stackrel{\mathcal{G}}{\rightarrow} & D^b(A) \\ \\
                      _{R_i \otimes_B -}\downarrow\;\;\;\;\;\;\;\;\;\;\;\;\;\;& & \;\;\;\;\;\;\;\;\;\;\;\;\;\;\downarrow \;_{\mathcal{G}(R_i) \otimes_A -} \\ \\
                      D^b(B) & \stackrel{\mathcal{G}}{\rightarrow} & D^b(A) \end{array} \]

This will give the image of $R_i$ in $TrPic(A)$ under the group isomorphism \[\Gamma : TrPic(B) \stackrel{\sim}{\rightarrow} TrPic(A)\] induced by $\mathcal{G}.$ \\ From now on, when we multiply two elements of $TrPic,$ we omit the tensor product  symbol $\otimes_A$ ($FF' := F \otimes_A F'$).

\begin{Proposition} \label{Gamma} We have
\[ \Gamma(R_i) = \left\{ \begin{array}{cc} F_i F_{i+1} F_i^{-1} &  \hbox{if} \;\; 1 \le i \le n-1, \\ F_n & \hbox{if} \;\; i = n. \end{array} \right. \]
\end{Proposition}

We need the following lemmas
\begin{Lemma} \label{need2} For all $1 \le i \le n-1,$
{\small \[ \mathcal{R}_i(T_j) = \left\{ \begin{array}{cccc}
               & T_j & \hbox{if} & |i-j| > 1\\
               &T_{i-1} & \hbox{if} & j = i-1 \\
             &T_{i+1} & \hbox{if} & j = i \\
 Cone(Cone(T_i  \rightarrow T_{i+1}) \rightarrow & T_{i+1}) & \hbox{if} & j = i+1 \end{array} \right. \]}
\end{Lemma}

\noindent {\bf Proof} For all $1 \le j \le n-1$ denote by
\[\begin{array}{ccccccccc}
\nu_{j,j+1} : &  Q_j & \rightarrow & Q_{j+1} &  \hbox{and} & \nu_{j+1,j} : & Q_{j+1} & \rightarrow & Q_j\\
              &  x   & \mapsto     & x(j \; j+1) &         &               & x       & \mapsto     & x(j+1 \; j) \end{array}\] and for all $1 \le i \le n$ denote by $\rho_i : Q_i \rightarrow Q_i$ where $\rho_i((i)) = (i).$ \\

For all $1 \le j \le n-1$ we have a triangle
\[ T_j \stackrel{f_j}{\rightarrow} T_{j+1} \stackrel{g_j}{\rightarrow} Q_j[j] \stackrel{h_j}{\rightsquigarrow} T_j[1] \]
where $f_j = (\rho_n, \rho_{n-1}, \cdots, \rho_{j+1}, 0, \cdots, 0),$ $\rho_k$ in degree $k-1;$ $g_j = \nu_{j+1,j}$ in degree $j$ and 0 otherwise; and $h_j = \rho_j$ in degree $j$ and 0 otherwise. \\

  Therefore, \[Cone(T_{j+1}[-1] \rightarrow Q_j[j-1]) \simeq T_j \;\; \hbox{and} \] \[Cone(T_j \rightarrow T_{j+1}) \simeq Q_j[j].\]
For $j > i + 1,$ we will prove the lemma using induction : if $j =
n,$
\[\mathcal{R}_i(T_n) = \mathcal{R}_i(Q_n[n-1]) = Q_n[n-1] = T_n.\]
Suppose that $\mathcal{R}_i(T_{j+1}) = T_{j+1}.$ We will show that
$\mathcal{R}_i(T_j) = T_j.$ {\small \begin{eqnarray}
\label{induction} \mathcal{R}_i(T_j) & \simeq &
\mathcal{R}_i(Cone(T_{j+1}[-1] \rightarrow Q_j[j-1])) \nonumber \\ &
\simeq & Cone(\mathcal{R}_i(T_{j+1})[-1] \rightarrow
\mathcal{R}_i(Q_j[j-1])) \nonumber \\ & \simeq & Cone(T_{j+1}[-1]
\rightarrow Q_j[j-1]) \simeq T_j. \end{eqnarray}} If $j = i +1,$
{\small \begin{eqnarray*} \mathcal{R}_i(T_{i+1}) & = &
\mathcal{R}_i(Cone(T_{i+2}[-1] \rightarrow Q_{i+1}[i])) \\ & \simeq
& Cone(\mathcal{R}_i(T_{i+2})[-1] \rightarrow
\mathcal{R}_i(Q_{i+1}[i])) \\ & \simeq & Cone(T_{i+2}[-1]
\rightarrow (Q_i \rightarrow Q_{i+1})[i]) \\ & \simeq & Cone \left(
\begin{array}{cccc} (Q_n \rightarrow \cdots \rightarrow & Q_{i+3} &
\stackrel{\nu_{i+3,i+2}}{\longrightarrow} & Q_{i+2}) \\ &
\downarrow\;_0 & & \;\;\;\;\;\;\;\;\;\;\downarrow \;_{\nu_{i+2,i+1}}
\\ & (Q_i & \stackrel{\nu_{i,i+1}}{\longrightarrow} & Q_{i+1})
\end{array} \right) \\ & & \hbox{where} \; Q_{i+1} \;\hbox{lies in
degree} \; i \\ & \simeq & \left(\begin{array}{ccccc} Q_n
\rightarrow \cdots \rightarrow Q_{i+3} & \rightarrow & Q_{i+2}& & \\
& \searrow \;_0 & \oplus & \searrow & \\ & &  Q_i &  \rightarrow &
Q_{i+1} \end{array}\right)  \\ & \simeq & Cone \left(
\begin{array}{cc} & Q_i \\  & \;\;\;\;\;\;\;\;\;\;\downarrow
\;_{\nu_{i,i+1}} \\  (Q_n \rightarrow \cdots \rightarrow Q_{i+2}
\rightarrow & Q_{i+1}) \end{array} \right) \\ & &  (\hbox{since} \;
Hom_B(Q_{i+3},Q_i) = 0) \\ & \simeq & Cone(Q_i[i] \rightarrow
T_{i+1}) \\ & \simeq & Cone(Cone(T_i  \rightarrow T_{i+1})
\rightarrow T_{i+1})  \end{eqnarray*}} If $j=i,$
{\small \begin{eqnarray*} \mathcal{R}_i(T_i) & = & \mathcal{R}_i(Cone(T_{i+1}[-1] \rightarrow Q_i[i-1])) \\ & \simeq & Cone(\mathcal{R}_i(T_{i+1}[-1]) \rightarrow \mathcal{R}_i(Q_i[i-1]))  \\ & \simeq & Cone \left( \left( \begin{array}{ccccc} Q_n \rightarrow \cdots \rightarrow Q_{i+3} &  \rightarrow  & Q_{i+2} & & \\ & \searrow & \oplus & \searrow & \\ & & Q_i &\rightarrow & Q_{i+1} \end{array} \right) \rightarrow Q_i[i] \right) \\ & &  \; \hbox{where} \; Q_{i+1} \;\hbox{lies in degree} \; i-1 \\
& \simeq & \left( \begin{array}{ccccc} Q_n \rightarrow \cdots \rightarrow Q_{i+3} &\stackrel{\nu_{i+3,i+2}}{\rightarrow}& Q_{i+2} &  & \\ & _0\searrow & \oplus & \;\;\;\;\;\;\;\;\;\;\;\;\searrow \;_{\nu_{i+2,i+1}} &  \\ & & Q_i & \stackrel{\nu_{i,i+1}}{\rightarrow}  & Q_{i+1} \\ & & & \searrow \;_{\rho_i} &  \oplus\\ & &  & & Q_i \end{array} \right)  \\ & &  \; \hbox{where} \; Q_{i+1} \;\hbox{lies in degree} \; i \\
& \simeq & Cone \left(  (Q_n \rightarrow \cdots \rightarrow Q_{i+3}
\rightarrow Q_{i+2}) \stackrel{\nu_{i+2,i+1}}{\longrightarrow}
\left( \begin{array}{ccc} Q_i &
\stackrel{\nu_{i,i+1}}{\longrightarrow} & Q_{i+1} \\  &
_{\rho_i}\searrow & \oplus \\  & & Q_i  \end{array}
\right)[i]\right) \\ & \simeq & Cone \left( T_{i+2}[-1] \rightarrow
Cone \left( \begin{array}{ccc} Q_i & \rightarrow & Q_{i+1} \\
\downarrow & & \\ Q_i & & \end{array} \right)[i-1] \right) \\ &
\simeq & Cone(T_{i+2}[-1] \rightarrow Q_{i+1}[i]) \simeq T_{i+1}.
\end{eqnarray*}} For $j < i,$ using induction and the same argument as
(\ref{induction}) we only need to show that
\[\mathcal{R}_i(T_{i-1}) = T_{i-1}.\] We have : {\small
\begin{eqnarray*} \mathcal{R}_i(T_{i-1}) & = &
\mathcal{R}_i(Cone(T_i[-1] \rightarrow Q_{i-1}[i-2])) \\ & \simeq &
Cone( \mathcal{R}_i(T_i[-1]) \rightarrow
\mathcal{R}_i(Q_{i-1}[i-2])) \\ & \simeq & Cone( T_{i+1}[-1]
\rightarrow ( Q_i \rightarrow Q_{i-1})[i-2]) \\ & \simeq & Cone
\left( \begin{array}{ccc} (Q_n \rightarrow \cdots \rightarrow &
Q_{i+1}) & \\ &\;\;\;\;\;\;\;\;\;\; \downarrow \; _{\nu_{i+1,i}}& \\
& (Q_i & \rightarrow Q_{i-1}) \end{array} \right)  \\ & &
\hbox{where} \; Q_{i-1} \;\hbox{in degree} \; i-2  \\ & \simeq &
T_{i-1}. \;\; \end{eqnarray*}} \qed

\begin{Lemma} \label{need3}
\[ \mathcal{R}_n(T_j) = \left\{ \begin{array}{cc} Cone(T_n \rightarrow T_j) & \hbox{if} \;\; j < n, \\ T_n[1] & \hbox{if} \;\; j= n. \end{array} \right. \]
\end{Lemma}

\noindent {\bf Proof} For $j=n,$ $R_n(T_n) = R_n(Q_n[n-1]) = Q_n[n]
= T_n[1].$ \\ For $j < n,$ we will show the result using induction.
{\small \begin{eqnarray*} \mathcal{R}_n(T_{n-1}) & = &
\mathcal{R}_n(Cone(Q_n \rightarrow Q_{n-1})) \\ & \simeq &
Cone(\mathcal{R}_n(Q_n) \rightarrow \mathcal{R}_n(Q_{n-1})) \\ &
\simeq & Cone \left( \begin{array}{ccc} Q_n & & \\ _{\rho_n}
\downarrow\;\;\;\;\;\; & & \\ (Q_n & \rightarrow & Q_{n-1})
\end{array} \right) \\ & \simeq & Cone(T_n \rightarrow T_{n-1})
\end{eqnarray*}} Assume that $\mathcal{R}_n(T_{j+1}) = Cone(T_n
\rightarrow T_{j+1}).$ We will show that \[\mathcal{R}_n(T_j) =
Cone(T_n \rightarrow T_j).\] {\small \begin{eqnarray*}
\mathcal{R}_n(T_j) & \simeq & \mathcal{R}_n(Cone(T_{j+1}[-1]
\rightarrow Q_j[j-1])) \\ & \simeq & Cone(\mathcal{R}_n(T_{j+1}[-1])
\rightarrow \mathcal{R}_n(Q_j[j-1])) \\ & \simeq & Cone(Cone(T_n
\rightarrow T_{j+1})[-1] \rightarrow Q_j[j-1]) \\ & \simeq &
Cone(T_n \rightarrow Cone(T_{j+1}[-1] \rightarrow Q_j[j-1])) \\ &
\simeq & Cone(T_n \rightarrow T_j). \;\; \end{eqnarray*} } \qed

\noindent {\bf Proof of Proposition \ref{Gamma}} We have a commutative diagram
\begin{eqnarray} \label{diagram} D^b(B) & \stackrel{\mathcal{G}}{\longrightarrow} & D^b(A) \nonumber \\ & &  \nonumber\\ \downarrow \; _{\mathcal{R}_i} & & \downarrow \; _{\mathcal{W}_i} \nonumber \\ & & \nonumber \\  D^b(B) & \stackrel{\mathcal{G}}{\longrightarrow} & D^b(A) \end{eqnarray}   where we denote  \[ \mathcal{W}_i := \mathcal{G} \circ \mathcal{R}_i \circ \mathcal{G}^{-1}. \]
We will show that for all  indecomposable projective $A$-modules $P,$
\begin{enumerate} \item $\mathcal{W}_i(P) = \mathcal{F}_i \mathcal{F}_{i+1} \mathcal{F}_i^{-1}(P)$ for all $1 \le i \le n-1,$ and
\item  $\mathcal{W}_n(P) = \mathcal{F}_n(P).$
\end{enumerate}
\begin{enumerate}
\item Fix $1 \le i \le n-1.$  Recall that $\mathcal{G}(T_j) = P_j$ for all $1 \le j \le n.$ By Lemma \ref{need2} and the above commutative diagram (\ref{diagram}),
{\small \begin{eqnarray*} \mathcal{W}_i(P_j) & = &  \mathcal{W}_i(\mathcal{G}(T_j)) \\ & = & \mathcal{G}\mathcal{R}_i(T_j) \\ & = &
\left\{ \begin{array}{ccccc}
& P_j & \hbox{if} & |i - j| > 1 \; \hbox{or}& j = i-1, \\
&  P_{i+1} & \hbox{if} & j = i, & \\  P_i \rightarrow P_{i+1} \rightarrow & P_{i+1} & \hbox{if} & j = i+1.& \end{array} \right.\end{eqnarray*}}
From (\ref{FiFjFi}) (see Theorem \ref{action}) we get
{\small \begin{equation} \label{W} \mathcal{F}_i \mathcal{F}_{i+1} \mathcal{F}_{i}^{-1}(P_j) =
\left\{ \begin{array}{ccccc}
& P_j & \hbox{if} & |i - j| > 1 \; \hbox{or}& j = i-1, \\
&  P_{i+1} & \hbox{if} & j = i, & \\  P_i \rightarrow P_{i+1} \rightarrow & P_{i+1} & \hbox{if} & j = i+1.& \end{array} \right.\end{equation}}
Therefore, we get $\mathcal{W}_i = \mathcal{F}_i \mathcal{F}_{i+1} \mathcal{F}_i^{-1}.$
\item By (\ref{GTi}), Lemma \ref{need3} and Diagram (\ref{diagram}),
{\small \begin{eqnarray*} \mathcal{W}_n(P_j) & = & \mathcal{W}_n\mathcal{G}(T_j) \\                                     & = & \mathcal{G}\mathcal{R}_n(T_j) \\
                                 & = & \mathcal{G}(Cone(T_n \rightarrow T_j)) \\ & = & P_n \rightarrow P_j \;\; \hbox{for all} \;\; 1 \le j < n, \end{eqnarray*}}and $\mathcal{W}_n(P_n) =  \mathcal{W}_n \mathcal{G}_n(T_n) = \mathcal{G} \mathcal{R}_n(T_n) = \mathcal{G}(T_n[1]) = P_n[1].$
Hence, by Lemma \ref{imageFi}, $\mathcal{W}_n = \mathcal{F}_n.$  \qed \end{enumerate}

\begin{Proposition} \label{gamma}
$\Gamma^{-1}(F_n) = R_n$ and
\begin{eqnarray*} \Gamma^{-1}(F_k) & = & R_k^{-1} R_{k+1}^{-1} \cdots R_{n-1}^{-1} R_n R_{n-1} \cdots R_{k+1}  R_k \\ & = & R_n R_{n-1} \cdots R_{k+1} R_k R_{k+1}^{-1} \cdots R_{n-1}^{-1} R_n^{-1} \end{eqnarray*} \end{Proposition}

\noindent {\bf Proof}
Define  $C_n := R_n$ and \[C_k := R_k^{-1} C_{k+1}  R_k, \;\; \hbox{for all} \;\; 1 \le k \le n-1.\]
By Proposition \ref{Gamma}, $\Gamma(C_n) = F_n.$ Now assume that $\Gamma(C_{k+1}) = F_{k+1}.$ We will show that $\Gamma(C_k) = F_k.$
\begin{eqnarray*} \Gamma(C_k) & = & \Gamma(R_k^{-1} C_{k+1} R_k) \\
                              & = & (F_k F_{k+1}^{-1} F_k^{-1}) F_{k+1} (F_k F_{k+1} F_k^{-1}) \\ & = & F_k F_{k+1}^{-1} F_k^{-1} F_k F_{k+1} F_k F_k^{-1} = F_k. \end{eqnarray*}  \qed

\begin{Corollary}The following diagram is commutative :
\[ \begin{array}{ccc} B(K_n) & \stackrel{\eta}{\longrightarrow} & B(A_n) \\ & & \\ \downarrow \; _{\varphi} & & \downarrow \; _{\psi} \\ & & \\ TrPic(A) & \stackrel{\Gamma^{-1}}{\longrightarrow} & TrPic(B) \end{array}\] This defines a categorification of the group homomorphism $\eta.$ \end{Corollary}

\begin{Corollary} Since $\eta$ is not injective, the action of $B(K_n)$ on $D^b(A)$ is not faithful. \end{Corollary}

\section{A faithful action}

Denote by $B(\tilde{A}_{n-1})$ the affine braid group generated by $h_1,\cdots,h_n$ subject to relations
\begin{eqnarray*} h_i h_{i+1} h_i & = & h_{i+1} h_i h_{i+1} \;\; \hbox{if} \; 1 \le i \le n-1,  \\ h_i h_j & = & h_j h_i \;\;\;\;\;\;\;\;\;\;\;\; \hbox{if} \; |i - j| \neq 1, n-1  \; \hbox{and} \\ h_n h_1 h_n & = & h_1 h_n h_1 . \end{eqnarray*}

In \cite{Schaps-Zakay-Illouz} Schaps and Zakay-Illouz obtained a group homomorphism between $B(\tilde{A}_{n-1})$ and the subgroup of $TrPic(A)$ generated by $H_1,\cdots,H_n$ where $H_i$ is the refolded tilting complex associated to the transposition $(i \;\; i+1).$ (See \cite{Schaps-Zakay-Illouz} for the definition).

\begin{Theorem} \label{SZI} (Schaps and Zakay-Illouz \cite{Schaps-Zakay-Illouz}) \\ There is a group homomorphism
\begin{eqnarray*} \rho :&  B(\tilde{A}_{n-1}) & \rightarrow  TrPic(A) \\
                           & h_i & \mapsto  H_i. \end{eqnarray*}
\end{Theorem}

{\bf Remark}  In \cite{Schaps-Zakay-Illouz} the action is defined for Brauer star algebras in general. \\

Denote by $\mathcal{H}_i$ the corresponding functor $\mathcal{H}_i := H_i \otimes_A -.$\\ Using \cite[Proposition 1]{Schaps-Zakay-Illouz} we get the images of indecomposable projective modules of $A$ under $\mathcal{H}_i$ :
\[\mathcal{H}_i(P_j) = \left\{ \begin{array}{ccccccc} & & & & P_j & \hbox{if}& |i - j| > 1, \\ & & & & P_{i-1} & \hbox{if}& j = i-1, \\& & & &  P_{i+1} & \hbox{if} &  j = i,  \\ P_{i}& \rightarrow &P_{i+1} &\rightarrow & P_{i+1} & \hbox{if} & j = i+1. \end{array} \right. \]
Comparing this with (\ref{W}) we get
\begin{equation} \label{h} \mathcal{H}_i = \mathcal{F}_i \mathcal{F}_{i+1} \mathcal{F}_i^{-1} \;\; \hbox{for all} \;\; 1 \le i \le n-1, \end{equation} and it is easy to see from (\ref{FiFjFi}) that
\begin{equation} \label{H} \mathcal{H}_n = \mathcal{F}_n \mathcal{F}_1 \mathcal{F}_n^{-1}. \end{equation}

\subsection{Injective braid group homomorphisms}
Let $B(B_n)$ be the braid group associated to Dynkin diagram of type $B_n$ with generators $b_1,b_2,\cdots,b_n$ subject to relations :
\begin{eqnarray*} b_i b_{i+1} b_i & = & b_{i+1} b_i b_{i+1} \;\; \hbox{if} \;\;  1 \le i \le n-1, \\ b_i b_j & = & b_j b_i \;\;\;\;\;\;\;\;\;\;\;\; \hbox{if}\;\; |i-j| > 1, \\ b_{n-1} b_n b_{n-1} b_n & = & b_n b_{n-1} b_n b_{n-1}. \end{eqnarray*}
 By \cite{Graham-Lehrer}, $B(B_n)$ is also generated by $\tau, s_1, s_2, \cdots, s_n$ subject to the relations :
\begin{eqnarray*} s_1 s_n s_1 & = & s_n s_1 s_n, \\ s_i s_{i+1} s_i & = & s_{i+1} s_i s_{i+1} \;\; \hbox{for} \;\; 1 \le i \le n-1, \\ s_i s_j & = & s_j s_i \;\;\;\;\;\; \hbox{for} \;\; |i-j| \neq 1,\;n-1, \\
                  \tau s_i \tau^{-1} & = & s_{i-1} \;\;\;\;\;\; \hbox{for} \;\; 1 \le i \le n-1, \\ \tau s_1 \tau^{-1} & = & s_n, \end{eqnarray*}
and the subgroup of $B(B_n)$ generated by $s_1, \cdots, s_n$ is the affine braid group $B(\tilde{A}_{n-1}).$ By \cite{Kent-Peifer}, this shows that $B(B_n)$ is a semidirect product of the infinite cyclic group generated by $\tau$ and $B(\tilde{A}_{n-1}).$ Therefore $B(\tilde{A}_{n-1})$ injects into $B(B_n)$ (see \cite{Kent-Peifer} for more details).

We can reformulate \cite[Proposition 6]{Graham-Lehrer} to get
\begin{eqnarray*} \tau & = & b_n b_{n-1} \cdots b_2, \\ s_i & = & b_i \;\; 1 \le i \le n-1,\\
                  s_n & = &  \tau b_1 \tau^{-1} =  b_n b_{n-1} \cdots b_2 b_1 b_2^{-1} \cdots b_{n-1}^{-1} b_n^{-1}, \end{eqnarray*} such that the embedding of $B(\tilde{A}_{n-1})$ into $B(B_n)$ is defined by the following map :
\begin{eqnarray*} \mu : & B(\tilde{A}_{n-1}) & \rightarrow B(B_n) \\
                         &      h_i     & \mapsto b_i \;\;\;\;\;\; 1 \le i \le n-1, \\
                         &      h_n  & \mapsto b_n b_{n-1} \cdots b_2 b_1 b_2^{-1} \cdots b_{n-1}^{-1} b_n^{-1}. \end{eqnarray*}

Now define a group homomorphism
\begin{eqnarray*} \chi : & B(B_n) & \rightarrow B(A_n) \\
                         &      b_i     & \mapsto \sigma_i \;\; 1 \le i \le n-1, \\
                         &      b_n  & \mapsto \sigma_n^{2}. \end{eqnarray*}
By \cite{Crisp},  $\chi$ is injective. Therefore we obtain an injective composition of injective group homomorphisms
\begin{eqnarray*} \chi \circ \mu : & B(\tilde{A}_{n-1}) & \rightarrow B(A_n) \\
                         &      h_i     & \mapsto \sigma_i, \;\; 1 \le i \le n-1, \\
                         &      h_n  & \mapsto \sigma_n  \sigma_n \sigma_{n-1} \cdots \sigma_2 \sigma_1 \sigma_2^{-1} \cdots \sigma_{n-1}^{-1} \sigma_n^{-1} \sigma_n^{-1}. \end{eqnarray*}

\subsection{Another Categorification}

By Proposition \ref{Gamma}, Proposition \ref{gamma}, (\ref{h}) and (\ref{H}) we obtain
\begin{eqnarray*} \Gamma^{-1}(H_i) & = & R_i \;\;\hbox{for all}\;\; 1 \le i \le n-1,\\
                  \Gamma^{-1}(H_n) & = & \Gamma^{-1}(F_n F_1 F_n^{-1}) \\
                                   & = & R_n  R_n R_{n-1} \cdots R_2 R_1 R_2^{-1} \cdots R_{n-1}^{-1} R_n^{-1} R_n^{-1}. \end{eqnarray*}

Therefore, we get the following commutative diagram :
\[ \begin{array}{ccc} B(\tilde{A}_{n-1}) & \stackrel{\chi \circ \mu}{\longrightarrow} & B(A_n) \\  \downarrow \;_{\rho} & & \downarrow \; _{\psi}  \\ TrPic(A) & \stackrel{\Gamma^{-1}}{\longrightarrow} & TrPic(B) \end{array} \]

And since $\chi \circ \mu$ and $\psi$ are injective, we obtain the following result :

\begin{Theorem}
The action of $B(\tilde{A}_{n-1})$ on $D^b(A)$ defined in Theorem \ref{SZI} is faithful.
\end{Theorem}

{\bf Remark} We leave as an open question whether the action defined in \cite{Schaps-Zakay-Illouz} is faithful for Brauer star algebras with exceptional vertex. We also get a faithful action of $B(B_n)$ on $D^b(A)$ by defining
\begin{eqnarray*} B(B_n) & \rightarrow & TrPic(A), \\
                    b_i  & \mapsto & F_i F_{i+1} F_{i} \;\; \hbox{if} \;\; 1 \le i \le n-1, \\  b_n     &  \mapsto & F_n F_n. \end{eqnarray*}

{\small {\sc Intan Muchtadi-Alamsyah \\ Department of Mathematics, \\University of Leicester, \\University Road, \\LE1 7RH, Leicester, UK,} \\{\it E-mail address} : idma1@le.ac.uk \\ \\ Current address :  \\{\sc Institut Teknologi Bandung,\\ Faculty of Mathematics and Natural Sciences, \\Algebra Research Group,\\ Jl.Ganesha no.10, \\Bandung 40132, Indonesia,} \\ {\it E-mail address} : ntan@math.itb.ac.id}

\end{document}